\documentclass[runningheads]{llncs}
\usepackage[T1]{fontenc}
\usepackage{graphicx}
\usepackage{algorithm}
\usepackage{algorithmicx}
\usepackage{algpseudocode}
\usepackage{array}
\usepackage{booktabs}
\usepackage{subcaption}
\usepackage{svg}
\usepackage{todonotes}
\usepackage[hidelinks]{hyperref}
\usepackage{tikz}
\usepackage{pgfplots} 
\usepackage{amssymb} 
\pgfkeys{/pgfplots/tuftelike/.style={
  semithick,
  tick style={major tick length=4pt,semithick,black},
  separate axis lines,
  axis x line*=bottom,
  axis line shift=4pt,
  xlabel shift=10pt,
  axis y line*=left,
  axis y line shift=10pt,
  }
  }

\newcommand{\repeatthanks}{\textsuperscript{\thefootnote}}
\begin{document}
\title{Efficient and Scalable Kernel Matrix
	Approximations using Hierarchical
	Decomposition}
\titlerunning{\(\mathsf{datafold}\) and \(\mathsf{GOFMM}\)} 
\author{}
\institute{
\url{}
\email{}}
%
%
%
\author{Keerthi Gaddameedi\thanks{Equal contribution, joint first author}\inst{1}\orcidID{0009-0003-0658-3884} \and
Severin Reiz\repeatthanks \inst{1}\orcidID{0000-0001-5752-4233} \and
Tobias Neckel \inst{1}\and
Hans-Joachim Bungartz\inst{1}\orcidID{0000-0002-0171-0712}}
\authorrunning{Gaddameedi*, Reiz*, Neckel, Bungartz}
%
\institute{Technical University of Munich, School of Computation, Information and Technology
\url{https://www.cs.cit.tum.de/sccs/}\\
\email{\{keerthi.gaddameedi,s.reiz\}@tum.de}}
\maketitle              

\begin{abstract}
With the emergence of Artificial Intelligence, numerical algorithms are moving towards more approximate approaches. For methods  such as PCA or diffusion maps, it is necessary to compute \textit{eigenvalues} of a large matrix, which may also be \textit{dense} depending on the kernel. 
A global method, i.e. a method that requires all data points simultaneously, scales with the data dimension \(N\) and not with the intrinsic dimension \(d\); the complexity for an exact dense eigendecomposition leads to $\mathcal{O}(N^{3})$.
We have combined the two frameworks, $\mathsf{datafold}$ and \(\mathsf{GOFMM}\). The first framework computes diffusion maps, where the \textit{computational bottleneck} is the eigendecomposition while with the second framework we compute the eigendecomposition \textit{approximately} within the iterative Lanczos method. A hierarchical approximation approach scales roughly with a runtime complexity of $\mathcal{O}(Nlog(N))$ vs. $\mathcal{O}(N^{3})$  for a classic approach.  
We evaluate the approach on two benchmark datasets -- scurve and MNIST -- with strong and weak scaling using OpenMP and MPI on \textit{dense} matrices with maximum size of \(100k\times100k\). 

\end{abstract}
\section{Introduction}
\subsection{Motivation}
Data-driven approaches to solve real-world problems have led to a rapid increase in data sizes. The potential of such approaches is limited by the current state of computational power. The memory requirements of dense matrices (i.e.~matrices with mostly non-zero entries) is \textit{$\mathcal{O}(N^{2})$}.  Similarly, the time complexity for operations such as \textit{mat-vec} is \textit{$\mathcal{O}(N^{2})$}. Therefore, these operations become computationally infeasible when the size of the matrices is large. As a solution to this, we aim to find low-rank approximations of these matrices using hierarchical algorithms. The fast multipole method ($\mathsf{GOFMM}$) is a novel algorithm for approximating dense symmetric positive definite (SPD) matrices so that the quadratic space and time complexity reduces to \textit{$\mathcal{O}(Nlog(N))$} with a small relative error. Dense SPD matrices appear in areas such as scientific computing, data analytics and statistical inference.
$\mathsf{GOFMM}$ is \textit{geometry-oblivious}, meaning that it does not require the geometry information or the knowledge of how the data has been generated \cite{1}. It just requires the distribution of data as input. 

Real-world problems require solving high dimensional data. Data-driven models assume an intrinsic geometry in the data, referred to as a manifold which can be used to extract essential information of lower dimension. $\mathsf{datafold}$ is a Python package that provides these data-driven models to find an explicit manifold parametrization for point cloud data \cite{5} by using kernel matrices. Kernels correspond to dot products in a high dimensional feature space, and one uses them to efficiently solve non-linear cases in machine learning \cite{3}.

In this paper we enable $\mathsf{datafold}$ functionalities to be used in conjunction with $\mathsf{GOFMM}$ to scale execution.

\subsection{Proposed approach} 
Manifold learning approaches learn the intrinsic geometry of high-dimensional data without the use of predetermined classifications (unsupervised learning). There are several manifold learning algorithms such as isomap\cite{isomap}, locally linear embedding\cite{locallylinearembed}, Hessian embedding\cite{hessianembed} etc., but we focus on \textit{diffusion maps}. Like PCA, diffusion maps also consists of a kernel matrix computation that describes the relation of data points in the space. A Markov chain is defined using the kernel matrix which is then decomposed to compute the eigenvalues and eigenvectors. These eigenvalues and eigenvectors are used to find a lower dimension than the dimension of the ambient space. 

$\mathsf{datafold}$ provides data-driven models based also on diffusion maps for finding a parametrization of manifolds in point cloud data and to identify non-linear dynamical systems from time series data \cite{5}. Since the eigendecomposition of the kernel matrix is very expensive, especially for huge matrices, hierarchical approaches are applied to be able to reduce the quadratic complexity to \textit{$\mathcal{O}(Nlog(N))$}.
 
The framework $\mathsf{GOFMM}$ provides hierarchical algorithms for large, dense, symmetric and positive-definite matrices. Let $K \in \mathbb{R}^{N\times N}$ be a dense kernel matrix for manifold data that is to be approximated. Let it also be symmetric and positive-definite. 
The goal is to find an approximation $\widetilde{K}$ such that the construction and any matrix-vector multiplications take only \textit{$\mathcal{O}(Nlog(N))$} work. The approximation must also satisfy the condition that the relative error between the approximated and exact matrix remains small, 
\begin{equation}\label{eq:1.1}
	\frac{||\widetilde{K} - K||}{||K||} \enspace \leq \enspace \epsilon, \quad 0 < \epsilon < 1,
\end{equation}
where $\epsilon$ is a user-defined tolerance. We use then the \textit{implicitly restarted Arnoldi} iteration to perform an iterative eigendecomposition. The matrix-vector multiplications in every iteration is performed using hierarchical methods from $\mathsf{GOFMM}$, where the dense matrix is compressed once at the beginning and then evaluated in each iteration. The relative error of the resulting eigenvalues with a reference solution is recorded. The scalability of the combined integrated software is tested on multiple cores.

\subsection{Related work and contributions}
There has been a growing interest on randomized computation of matrix decompositions~\cite{MARTINSSON201147,doi:10.1137/100804139}. They also occur in theoretical deep learning, for example, with shallow Gaussian processes~\cite{garriga2018deep} or for finding weights in deep neural networks by solving a system of linear equations~\cite{bolager2023sampling} or for Hessian approximations in second-order optimization~\cite{reiz2022neural}. Naturally, approximate matrix calculations are suitable for data aplications, especially when the \textit{modelling} error (e.g., of neural networks) are bigger than the \textit{numerical} error. However, often matrices like kernels from radial basis functions, may not be global low-rank and only allow for low-rank treatment for off-diagonal matrices with the so-called \(\mathcal{H}\)-arithmetic~\cite{11,grasedyck2013literature}. Hence, our target here is matrices that have globally significant rank, but allow for approximations on the off-diagonals. To our knowledge, the most prominent framework for hierarchical structured matrices is STRUMPACK~\cite{doi:10.1137/20M1349667}; $\mathsf{GOFMM}$~\cite{mpigofmm} shows some superiority for kernel matrices against STRUMPACK, underlining that $\mathsf{GOFMM}$ a good candidate for diffusion maps kernels. In addition, eigendecompositions of dense kernel matrices are the computational bottleneck of diffusion maps, limiting the global size. Existing work from our group integrated the $\mathsf{GOFMM}$  and the $\mathsf{datafold}$  frameworks. 

\textbf{Contributions} of this paper include  
\begin{enumerate}
    \item To our knowledge, first \(\mathcal{H}\)-arithmetic in iterative eigendecompositions
    \item Analysis of dense kernel matrices from diffusion maps \textit{enabling} bigger sizes
    \item Versatile approach in software engineering to allow for better reproducibility and portability
\end{enumerate}
Our approach is \textit{using} the framework $\mathsf{GOFMM}$, and \textit{extends} $\mathsf{datafold}$ by offering a hierarchical variant for the eigendecomposition.


\section{Methods}
\subsection{Diffusion maps}\label{sec:diffusionmaps}
Diffusion Maps is a non-linear technique of dimensionality reduction. It tries to obtain
information about the manifold encoded in the data without any assumptions on the
underlying geometry. As opposed to using the Euclidean distance or the geodesic distance in isomaps, diffusion maps use an affinity or similarity matrix obtained by  using a kernel function
that produces positive and symmetric values. Given a dataset $X = \{x_1, x_2, x_3,\ldots ,x_n\}$ and a Gaussian kernel function, a similarity matrix can be computed as
\begin{equation}
W_{ij} \enspace = \enspace w(i,j) \enspace = e^{\frac{-||x_{i}-x_{j}||_{2}^{2}}{\sigma^{2}}},
\end{equation}
where $x_i$ and $x_j$ are a pair of data points and $\sigma$ is the radius of the neighborhood around the point $x_i$. As outlined in Algorithm \ref{alg:dmaps}, the similarity matrix is then normalized with the density $Q$ (degree of vertex) and the  density parameter $\alpha$ to capture the influence of the data distribution on our approximations. For $\alpha = 0$, the density has maximal influence on how the underlying geometry is captured and vice versa for $\alpha = 1$. Therefore, normalization is done with $\alpha = 1$, and then a Markov chain is defined to obtain the probabilities of transitioning from one point to another.
\begin{algorithm}
    \caption{DiffusionMaps \cite{COIFMAN20065}}\label{alg:dmaps}
    \begin{algorithmic}[1]
        \State Compute $W_{ij}$                   \Comment{Similarity matrix}
        \State Compute normalized weights $W_{ij}^{\alpha} \enspace = \enspace\frac{W_{ij}}{Q_{i}^{\alpha} \cdot Q_{j}^{\alpha}}$      \Comment{$Q^\alpha$: Influence of density}
        \State Define Markov chain $P_{ij} \enspace = \enspace \frac{W_{ij}^{\alpha}}{Q_{i}^{\alpha}}$                \Comment{P: Transition matrix}
        \State Perform $t$ random walks to obtain $P^t$      
        \State Perform eigendecomposition on $P^t$     \Comment{$\lambda_{r}$: eigenvalues, $\psi_{r}$: eigenvectors}
        \State Lower dimension $d(t)$ = max\{ $l$ : $\lambda_{l}^{t} < \delta \lambda_{1}^{t}$ \}     \Comment{$\delta$: Predetermined precision factor} 
    \end{algorithmic}
\end{algorithm}
 Then, the transition matrix $P^t$ is obtained by performing random walks for $t$ time steps. Afterwards, an eigendecomposition is performed on the transition matrix to compute the eigenpairs which are further used to obtain the underlying lower dimension of the dataset. The computational complexity of diffusion algorithms  in standard form is $O(N^3)$, and the  eigendecomposition is the most expensive part  of the algorithm. Hence, we tackle this by using hierarchical matrix approximations.
\subsubsection{Hierarchical partitioning}
If $K$ is a kernel matrix, the hierarchically low-rank approximation $\widetilde{K}$ of $K$ is given as \cite{10,11}
\begin{equation}\label{eq:2.18}
	\widetilde{K} \enspace = \enspace D \enspace + \enspace S \enspace + \enspace UV,
\end{equation} 
where $D$ is a block-diagonal matrix with every block being an hierarchical matrix (short: $\mathcal{H}$-matrix), $S$ is a sparse matrix and $U$, $V$ are low rank matrices. The $\mathcal{H}$-matrix $\widetilde{K}$ is to be computed such that the error from Equation \ref{eq:1.1} ranges in the order of the user defined tolerance $0 < \epsilon < 1$.
The construction of $\widetilde{K}$ and matrix-vector product both take $\mathcal{O}(N\log N)$ operations. We then incorporate these hierarchical approximations into the diffusion maps algorithm to improve the computational costs of the eigendecompositions. 
\subsection{Implicitly restarted Arnoldi method}
Implicit restarted Arnoldi method is a variation of Arnoldi process which builds on the \textit{power iteration} method which computes $Ax, Ax^{2},Ax^{3}...$ for an arbitrary vector $x$, until it converges to the eigenvector of the largest eigenvalue of matrix $A$. To overcome the drawbacks of so many unnecessary computations for a single eigenvalue and its corresponding eigenvector, the Arnoldi method aims to save the successive vectors as they contain considerable information that can be further exploited to find new eigenvectors. The saved vectors form a $Krylov$ matrix which is given as \cite{19}
\begin{equation}
	\mathcal{K}_{n} = \mathsf{Span}[ x, Ax, A^{2}x ... A^{n-1}x ].
\end{equation}
Orthonormal vectors $x_{1}, x_{2}, x_{3}...$ that span a \textit{Krylov} subspace are extracted using \textit{Gram-Schmidt} orthogonalization from each column of $Krylov$ matrix. The $k$-step Arnoldi iteration is given in Algorithm \ref{alg:arnoldi} \cite{18}.
\begin{algorithm}
	\caption{$k$-step ArnoldiFactorization($A$,$x$)}\label{alg:arnoldi}
	\begin{algorithmic}[1]
		\State $x_{1} \leftarrow \frac{x}{||x||}$\Comment{Computes first Krylov vector $x_1$}
		\State $w \leftarrow Ax_{1}$ \Comment{Computes new candidate vector}
		\State $\alpha_{1} \leftarrow x_{1}^{H}w$   
		\State $r_{1} \leftarrow w - \alpha_{1} x_{1}$ 
		\State $X_{1} \leftarrow [x_{1}]$ 				\Comment{Orthonormal basis of Krylov subspace}
		\State $H_{1} \leftarrow [\alpha_{1}]$  		\Comment{Upper Hessenberg matrix}						
		\ForAll{$ j = {1...k-1}$}     \Comment{For $k$ steps, compute orthonormal basis $X$}
		\State \Comment{and the projection of matrix $A$ on the new basis}
		\State $\beta_{j} \leftarrow ||r_{j}||$ ;\enspace $x_{j+1} \leftarrow \frac{r_{j}}{\beta_{j}}$ 
		\State $X_{j+1} \leftarrow [X_{j}, x_{j+1}]$ ; \enspace $\hat{H}_{j} \leftarrow \Big[H_{j}, \enspace \beta_{j} e_{j}^{T}\Big]^{T} $
         \State \Comment{$e_j$ is the standard basis of  coordinate vector space}
		\State $z \leftarrow Ax_{j}$ 
		\State $h \leftarrow X_{j+1}^{H}z$; \enspace $r_{j+1} \leftarrow z - X_{j+1}h$ \Comment{Gram-Schmidt Orthogonalization}
		\State $H_{j+1} \leftarrow [\hat{H}_{j}, h]$
		\EndFor
	\end{algorithmic}
\end{algorithm}
$H$ is the orthogonal projection of $A$ in the \textit{Krylov} subspace. It is observed that eigenvalues of the upper Hessenberg matrix $H$ (the so-called \textit{Ritz} values) converge to the eigenvalues of $A$. When the current iterate $r_{j}$ = 0, the corresponding Ritz pair becomes the eigenpair of $A$.\\\\
One of the drawbacks of Arnoldi process is that the number of iterations taken for convergence is not known prior to the computation of well-approximated Ritz values \cite{18}. This causes the computation of the Hessenberg matrix to be of complexity $\mathcal{O}(k^{3})$ at the $k$-th step. 
A more efficient approach is \textit{implicitly restarted Arnoldi method} uses an implicitly shifted \textit{QR-iteration}. It avoids storage and numerical instabilities associated with the standard approach by compressing the necessary information from very large \textit{Krylov} subspace into a fixed size $k$-dimensional subspace. \\\\
The Arnoldi factorization of length $m=k+p$ has the form
\begin{equation}\label{eq:arnoldi}
	AX_{m} = X_{m}H_{m} + r_{m}e_{m}^{T} \ .
\end{equation}
The implicit restarting method aims to compress this to length $k$ by using QR steps to apply $p$ shifts resulting in \cite{18}
\begin{equation}
	AX_{m}^{+} = X_{m}^{+}H_{m}^{+} + r_{m}e_{m}^{T}Q \ ,
\end{equation}
where $V_{m}^{+} = V_{m}Q$, $H_{m}^{+} = Q^{T}H_{m}Q$ and $Q=Q_{1}Q_{2}...Q_{p}$. $Q_j$ is the orthogonal matrix associated with the corresponding shift $\mu_{j}$. The first $k-1$ values of $e_{m}^{T}Q$ are zero and thus the factorization becomes
\begin{equation}
	AX_{k}^{+} = X_{k}^{+}H_{k}^{+} + r_{k}^{+}e_{k}^{T} \ .
\end{equation}
The residual $r_{m}^{+}$ can be used to apply $p$ steps to obtain back the $m$-step form. A polynomial of degree $p$ of the form $\prod_{1}^{p}(\lambda - \mu_{j})$ is obtained from these shifts. The roots of this polynomial are used in the QR process to filter components from the starting vector.
\subsubsection{Implicitly restarted Lanczos method}
Consider the Equation (\ref{eq:arnoldi}) for Arnoldi factorization. $X_{m}$ are orthonormal columns and $H_{m}$ is the upper Hessenberg matrix. If $A$ is a Hermitian matrix, it becomes Lanczos factorization. So Arnoldi is basically a generalization to non-hermitian matrices. For Lanczos method, $H_{m}$ is a real, symmetric and tridiagonal matrix and the $X_{m}$ are called Lanczos vectors. The algorithms hence remain the same as the ones described for Arnoldi. The method $\mathsf{scipy.sparse.linalg.eigs}$ uses Arnoldi iteration since it deals with real and symmetric matrices while $\mathsf{scipy.sparse.linalg.eigsh}$ invokes implementation of Lanczos methods.

\section{Implementation}
Manifold learning data is generated in Python using $\mathsf{datafold}$ and then the diffusion maps algorithm is invoked. The eigendecompositions contained in diffusion maps are performed using an subclass of $\mathsf{LinearOperator}$. $\mathsf{LinearOperator}$ is instantiated with a $\mathsf{matvec}$ implementation from $\mathsf{GOFMM}$. This is done by writing an interface using the Simplified Wrapper Interface Generator (SWIG\footnote{\url{swig.org}}) to access the $\mathsf{GOFMM}$ methods written in C++ from a Python script. In this section, we further delve into the details of how each part has been implemented. 
\subsection{Integration of $\mathsf{datafold}$ and \(\mathsf{GOFMM}\)} 
 The software architecture of $\mathsf{datafold}$ contains integrated models that have been implemented in a modularized fashion and an API that has been templated from $\mathsf{scikit-learn}$ library. The architecture as shown in Figure \ref{fig:workflow} consists of three layers and describes the hierarchy of the workflow.
\begin{figure}[htpb]
	\centering 
	\begin{tikzpicture}[node distance=4cm]
		\node[draw,
		minimum width=2.5cm,
		minimum height=1.5cm,text width=3cm] (block1) {layer 1 \\ \textbf{datafold.appfold}};
		
		\node[draw,
		right of=block1,
		minimum width=2.5cm,
		minimum height=1.5cm,text width=3cm] (block2) {layer 2 \\ \textbf{datafold.dynfold}};
		
		\node[draw,
		right of=block2,
		minimum width=2.5cm,
		minimum height=1.5cm,text width=3cm] (block3) {layer 3 \\ \textbf{datafold.pcfold}};
		\draw[-latex] (block1) edge (block2);
		\draw[-latex] (block2) edge (block3);
	\end{tikzpicture}
	\caption[Workflow Hierarchy of \textit{$\mathsf{datafold}$}]{Workflow hierarchy of $\mathsf{datafold}$}\label{fig:workflow}
\end{figure}
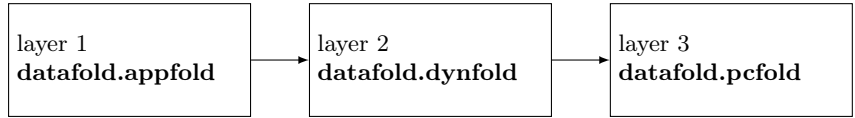\\
$\mathsf{datafold.appfold}$ is the highest level in the workflow hierarchy and contains meta-models that provide access to multiple sub-models captured in the class. The second layer $\mathsf{datafold.dynfold}$ provides models that deal with point cloud manifold or the dynamics of time series data. Finally, the last layer $\mathsf{datafold.pcfold}$ consists of fundamental algorithms such as eigensolvers, distance matrix computations etc. along with objects and data structures associated with them. The software maintains a high degree of modularity with this workflow and therefore allows usage of each layer's methods to be used on their own.
We have a docker file with commands to install all the run-time dependencies followed by installation of \(\mathsf{GOFMM}\) and $\mathsf{datafold}$. 
The docker image containing $\mathsf{GOFMM}$ and $\mathsf{datafold}$ has been converted to $\mathsf{Charliecloud}$ in order to be viable with the linux cluster at $\mathsf{LRZ}$. The docker image is then converted to a charliecloud image using the command \texttt{ch-builder2tar <docker-image> /dir/to/save}. Then, the charliecloud image is exported to the linux cluster and unpacked with the command \texttt{ch-tar2dir <charliecloud-image> /dir/to/\-unpack}. Once the compressed image is unpacked, the environment variables are set and \(\mathsf{GOFMM}\) is compiled. Finally, the SWIG interface file is compiled to generate Python versions of \(\mathsf{GOFMM}\)’s C++ methods. \\
\subsection{LinearOperator}\label{sec:3.3}
$\mathsf{SciPy}$ \cite{16} is an open-source free software with modules for common tasks of scientific computing such as linear algebra, solvers, interpolation etc. 
It contains seven matrix and array classes for different types of representations such as sparse row matrix, column matrix, coordinate format etc. It also accommodates methods to build various kinds of sparse matrices and two submodules $\mathsf{csgraph}$ and $\mathsf{linalg}$. The submodule $\mathsf{linalg}$ provides an abstract interface named $\mathsf{LinearOperator}$ that uses iterative solvers to perform matrix vector products. This interface consists of methods such as $\mathsf{matmat(x)}$, $\mathsf{matvec(x)}$, $\mathsf{transpose(x)}$ for matrix-matrix multiplication, matrix-vector multiplication and transposition of a matrix. A concrete subclass of $\mathsf{LinearOperator}$ can be built by implementing either one of $\mathsf{\_matvec}$ or $\mathsf{\_matmat}$ methods and the properties $\mathsf{shape}$ and $\mathsf{dtype}$. Depending on the type of matrices at hand, corresponding matvec methods may also be implemented.\\
\begin{table}[htpb]
	\caption[Matrix Factorization methods in $\mathsf{scipy.sparse.linalg}$]{Matrix Factorizations in $\mathsf{scipy.sparse.linalg}$.}\label{tab:factorizations}
	\centering
	\begin{tabular}{|l || m{7.9cm}|}
		\toprule
		$\mathsf{scipy.sparse.linalg}$$\mathsf{\textbf{.eigs}}$ & Computes eigenvalues and vectors of square matrix  \\
		\midrule
		$\mathsf{scipy.sparse.linalg}$$\mathsf{\textbf{.eigsh}}$ & Computes eigenvalues and vectors of real symmetric or complex Hermitian matrix  \\
		\midrule
		$\mathsf{scipy.sparse.linalg}$$\mathsf{\textbf{.lobpcg}}$ & Locally Optimal Block Preconditioned Conjugate Gradient Method  \\
		\midrule
		$\mathsf{scipy.sparse.linalg}$$\mathsf{\textbf{.svds}}$ & Partial Singular Value Decompositions  \\	
		\midrule
		$\mathsf{scipy.sparse.linalg}$$\mathsf{\textbf{.splu}}$ & LU decomposition of sparse square matrix  \\	
		\midrule
		$\mathsf{scipy.sparse.linalg}$$\mathsf{\textbf{.spilu}}$ & Incomplete LU decomposition of sparse square matrix  \\
		\midrule
		$\mathsf{scipy.sparse.linalg}$$\mathsf{\textbf{.SuperLU}}$ & LU decomposition of a sparse matrix  \\		
		\bottomrule
	\end{tabular}
\end{table}
$\mathsf{scipy.sparse.linalg}$ also provides methods for computing matrix inverses, norms, decompositions and linear system solvers. The functionality we are interested in are the matrix decompositions. In Table \ref{tab:factorizations}, we can take a look at various decomposition methods that are present in the module. The method we use to decompose data obtained from $\mathsf{datafold}$ is $\mathsf{scipy.sparse.linalg.eigsh}$ \cite{16}. This methods requires either an $\mathsf{ndarray}$, a sparse matrix or $\mathsf{LinearOperator}$ as parameters. It optionally takes $k$, which is the number of desired eigenvalues and eigenvectors. It solves $Ax[i] = \lambda_{i} x[i]$ and returns two arrays - $\lambda_i$ for eigenvalues and k vectors $X[:i]$, where $i$ is the column index corresponding to the eigenvalue. 

$\mathsf{scipy.sparse.linalg.eigsh}$ is a wrapper for the ARPACK functions SSEUPD and DSEUPD which use the implicitly restarted Lanczos method to solve the system for eigenvalues and vectors \cite{17}.

\section{Results}

Several experiments have been performed using datasets such as uniform distribution\footnote{\url{https://numpy.org/doc/stable/reference/random/generated/numpy.random.uniform.html}}, s-curve\footnote{\url{https://scikit-learn.org/stable/modules/generated/sklearn.datasets.make\_s\_curve.html} }, swiss-roll\footnote{\url{https://scikit-learn.org/stable/modules/generated/sklearn.datasets.make\_swiss\_roll.html}} and MNIST\cite{deng2012mnist}. 
Accuracy measurements for the datasets s-curve and MNIST have been presented in \autoref{sec:eigenpair}. 
Accuracy has been measured by computing Frobenius norm between eigenvalue computations of $\mathsf{scipy}$ solver and $\mathsf{GOFMM}$ and additionally, resultant eigenvectors have been plotted to provide a qualitative analysis. 
Furthermore, experiments were conducted to analyze performance through both weak and strong scaling in \autoref{sec:scaling}. Due to varying computational requirements, weak scaling experiments have been conducted on CoolMUC-2 linux cluster of LRZ\footnote{\url{https://doku.lrz.de/coolmuc-2-11484376.html}} and strong scaling on the supercomputer SuperMUC-NG\footnote{\url{https://doku.lrz.de/hardware-of-supermuc-ng-11482553.html}}. Efficiency and scalability of our approach were analyzed by examining results obtained with large problem sizes.
\subsection{Eigenvalue and eigenvector computations}\label{sec:eigenpair}
The experiments were performed on the CoolMUC-2 cluster of the Leibniz Supercomputing Centre\footnotemark[5] . It has 812 28-way Intel Xeon E5-2690 v3 ("Haswell") based nodes with 64GB memory per node and FDR14 Infiniband interconnect.
\subsubsection{S-curve}
A 3D S-curve dataset\footnotemark[3] is generated using $\mathsf{scikit-learn}$ \cite{scikit-learn} with 16384 points in the dataset. A 3D S-curve has an underlying intrinsic dimension of 2 and we apply diffusion maps algorithm to compute this. Since our focus lies in the eigendecompositions of the kernel matrix, eigenpairs are computed using two solvers. The first set of values are computed using the $\mathsf{scipy}$ solver and these are taken as reference values. The approximations of our  $\mathsf{GOFMM}$ $\mathsf{matvec}$ implementation are computed, and the error values in the Frobenius norm are observed to be in the range of \(9e-4\). 
\begin{figure}[htpb]
	\centering
	\begin{subfigure}{.5\textwidth}
		\centering
		\includegraphics[width=1.\linewidth]{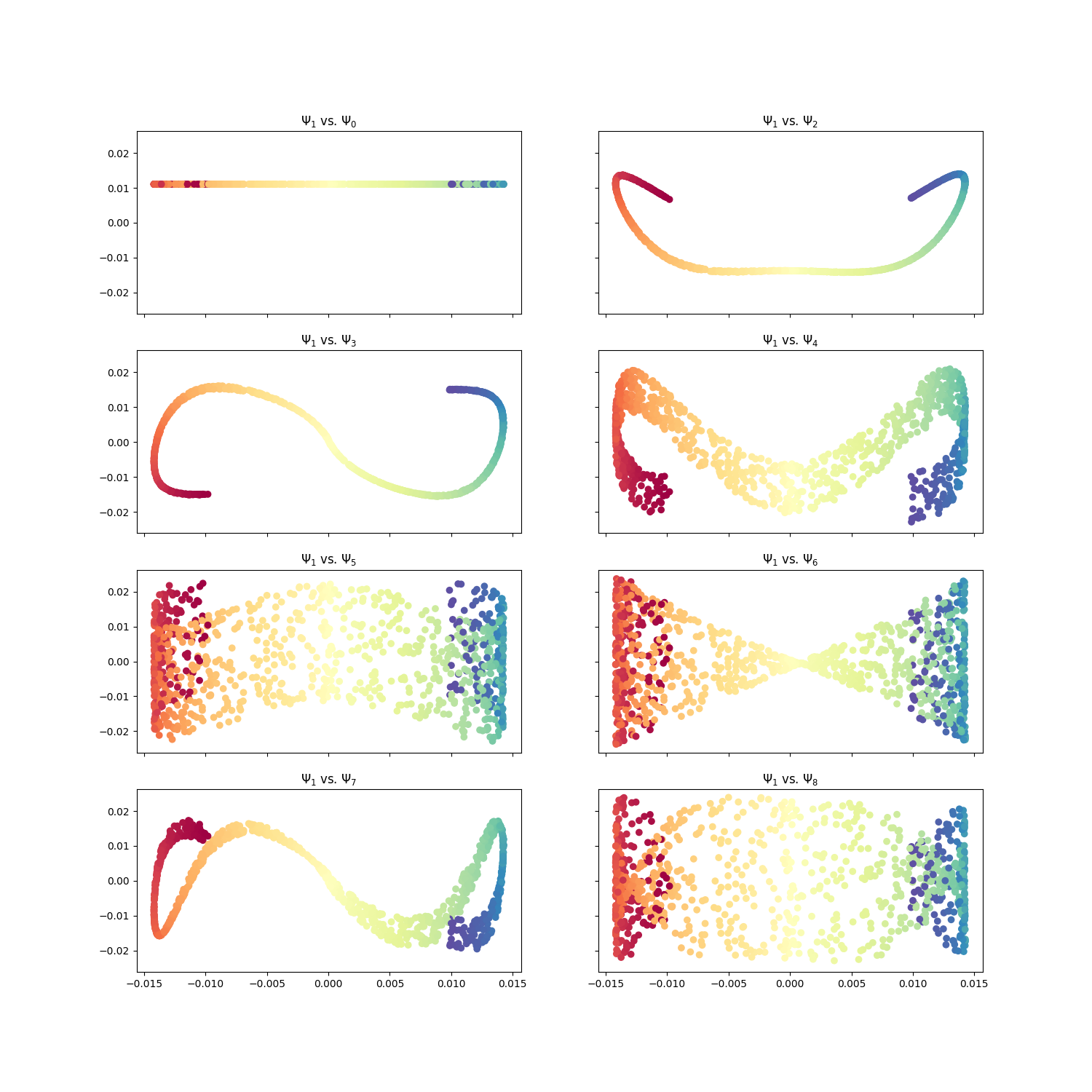}
	\end{subfigure}%
	\begin{subfigure}{.5\textwidth}
		\centering
		\includegraphics[width=1.\linewidth]{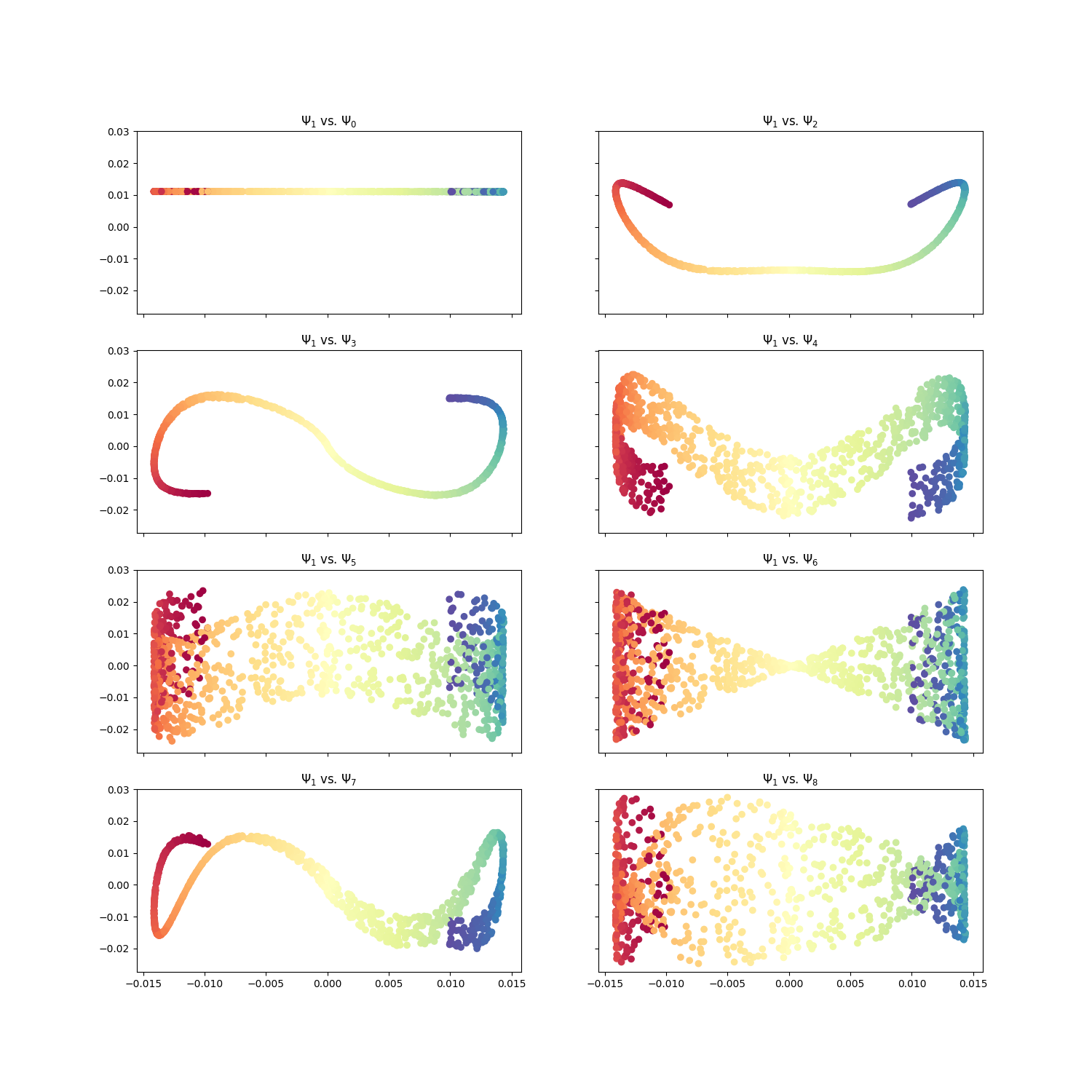}
	\end{subfigure}
	\caption{Eigenvector comparison for $\mathsf{scipy}$ solver on the left and $\mathsf{GOFMM}$ on the right for \textbf{scurve} }
	\label{fig:scurve}
\end{figure}
We can compare the embeddings obtained from both solvers by fixing the first non-trivial eigenvector and comparing it to the other eigenvectors. Eigenvector comparison for both $\mathsf{scipy}$ solver and $\mathsf{GOFMM}$ can be observed to be very similar in Figure \ref{fig:scurve}.
\subsubsection{MNIST}
\begin{figure}[H]
	\centering
	\begin{subfigure}{.5\textwidth}
		\centering
		\includegraphics[width=1.\linewidth]{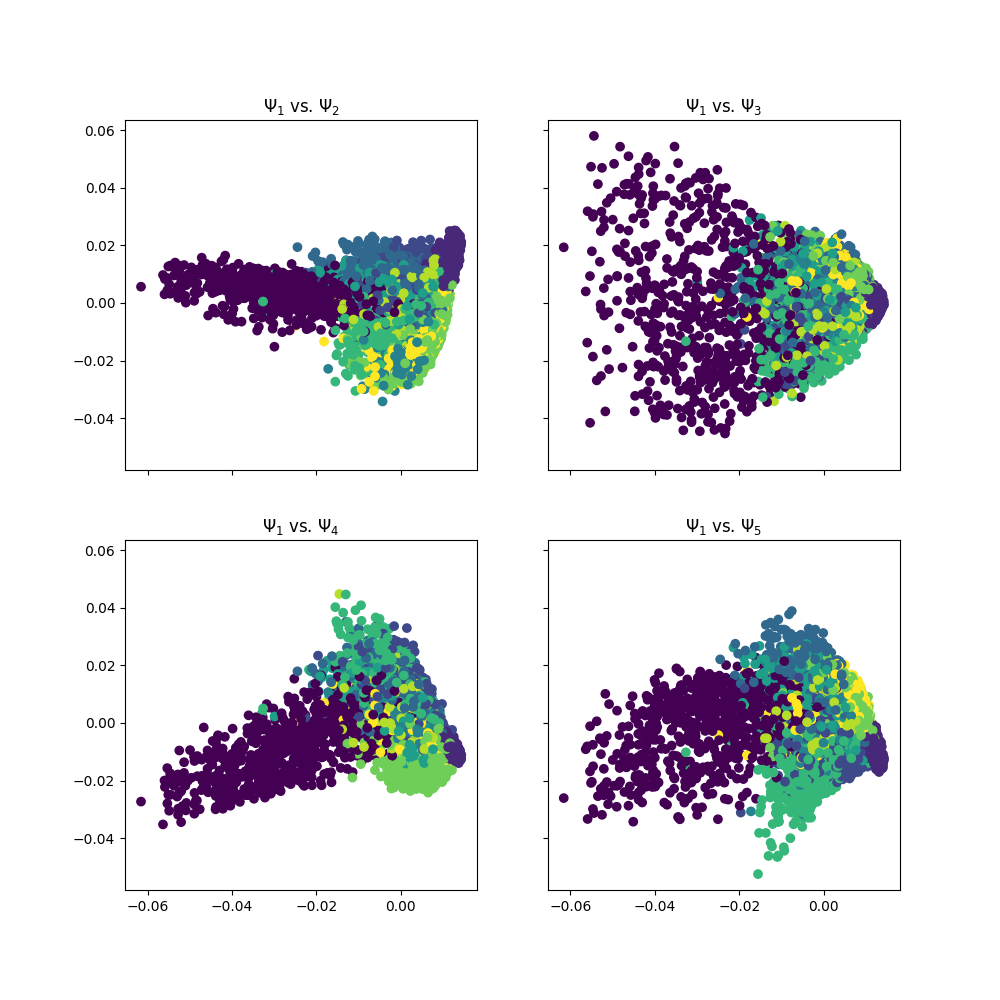}
	\end{subfigure}%
	\begin{subfigure}{.5\textwidth}
		\centering
		\includegraphics[width=1.\linewidth]{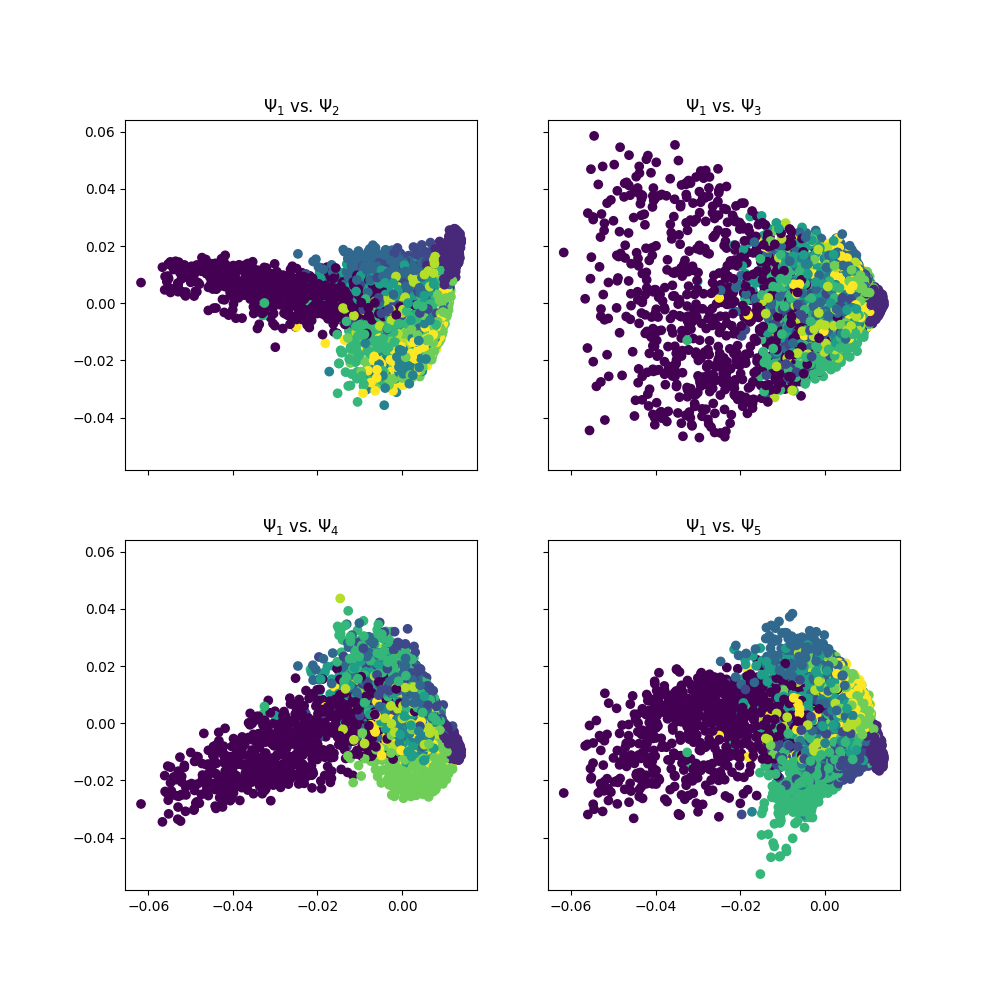}
	\end{subfigure}
	\caption{Eigenvector comparison for scipy solver on the left and \(\mathsf{GOFMM}\) on the right}
	\label{fig:hr_digits}
\end{figure}

The MNIST database (Modified National Institute of Standards and Technology database) \cite{deng2012mnist} is a large database of handwritten digits that is commonly used for training various image processing systems. MNIST has a testing sample size of 10,000 and a training size of 60,000 where each sample has 784 dimensions.\\
Due to a large dataset with 784 dimensions for each sample, MNIST makes a fitting application for hierarchical algorithms. Sample sizes of up to $16384$ are loaded from MNIST followed by diffusion maps algorithm applied to the dataset resulting in a kernel matrix of size $16k\times16k$. As previously mentioned, the goal is to perform efficient eigendecompositions using hierarchical algorithms. Therefore, eigenpairs are computed using $\mathsf{scipy}$ and $\mathsf{GOFMM}$ and we observe that for a matrix size of 8192, eigenvector comparison for both solvers look qualitatively similar as can be seen in Figure \ref{fig:hr_digits}. The Frobenius norm of the difference of the first five eigenvalues is also in the range of \(1e-4\). 
The parameters required to obtain the results show that the approach is very problem-dependent. As already mentioned in \cite{1}, problems with dense matrices are better suited to hierarchical approaches.

\subsection{Scaling}\label{sec:scaling}
\subsubsection{Complexity analysis}
As we have established previously in \ref{sec:diffusionmaps}, computational bottleneck of diffusion maps algorithm (e.g. for manifold learning, see Algorithm \ref{alg:dmaps}) is the \textit{eigenvector} (EV) computations. In general, for a matrix of size \(N\times N\), EV computations scale with a complexity of \(\mathcal{O}(N^3)\). In the past, matrices with large sizes in $\mathsf{datafold}$ were restricted to sparse matrices. A sparse matrix only requires \(\mathcal{O}(N)\) operations per iteration as one assumes a constant number of non-zero entries. Therefore with N rows, \textit{sparse matrix-vector multiplication} operation only costs \(\mathcal{O}(N)\). 
We usually also limit the number of iterations necessary for the Arnoldi method to a factor of  {\raise.17ex\hbox{$\scriptstyle\sim$}}100, resulting in an overall computational complexity of \(\mathcal{O}(N)\). 

However, there exist numerous kernels that do not result in sparse matrices and hence dense matrices are necessary. For a hierarchical approximate \textit{dense matrix-vector multiplication}, we need around \(\mathcal{O}(N \log N)\) operations. 

 Including FLOP counts in \cite{1,23}, we have looked at performance measurements for problem sizes up to 200k\footnote{ $\mathsf{GOFMM}$ can work with dense
matrices of 200k starting with at least 2 nodes, prohibiting the use of a bigger
matrix size for strong scaling analysis.}. Owing to the need for a high number of nodes, scaling experiments were performed on the Intel Xeon Platinum 8174 ("Skylake")  partition of SuperMUC-NG\footnotemark[6] which has 6,336 thin nodes and 144 fat nodes with 48 cores per node.  
\subsubsection{Weak scaling}
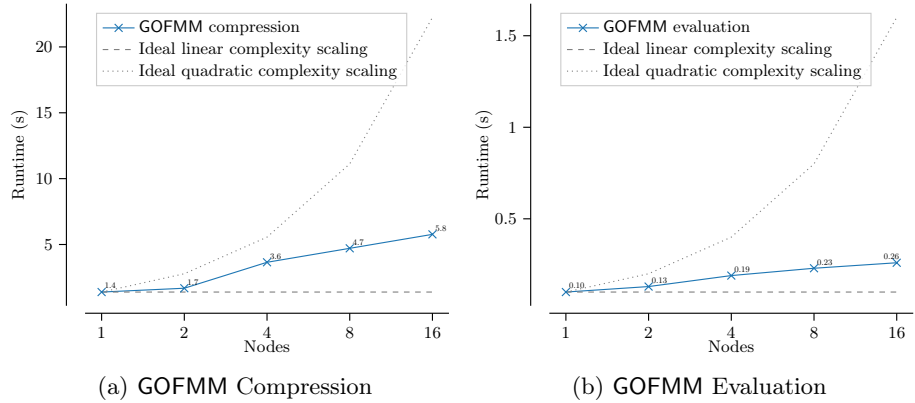
\begin{figure}
	\begin{subfigure}{.5\textwidth}
\resizebox{0.98                                                                                                 \textwidth}{!}{
\begin{tikzpicture}

\definecolor{color0}{rgb}{0.12156862745098,0.466666666666667,0.705882352941177}

\begin{axis}[tuftelike,
legend cell align={left},
legend style={
  fill opacity=0.8,
  draw opacity=1,
  text opacity=1,
  at={(0.03,0.97)},
  anchor=north west,
  draw=white!80!black
},
log basis x={10},
tick align=outside,
tick pos=both,
x grid style={white!69.0196078431373!black},
xlabel={Nodes},
xmin=0.870550563296124, xmax=18.3791736799526,
xmode=log,
xtick style={color=black},
xtick={1,2,4,8,16},
xticklabels={1,2,4,8,16}, 
y grid style={white!69.0196078431373!black},
ylabel={Runtime (s)},
ymin=0.3475, ymax=23.2825,
ytick style={color=black},
axis y line*=left,  
]
\addplot [semithick, color0, mark=x, mark size=3, mark options={solid}]
table {%
1 1.39
2 1.68
4 3.65
8 4.71
16 5.77
};
\addlegendentry{$\mathsf{GOFMM}$ compression}
\addplot [semithick, white!50.1960784313725!black, dashed]
table {%
1 1.39
2 1.39
4 1.39
8 1.39
16 1.39
};
\addlegendentry{Ideal linear complexity scaling}
\addplot [semithick, white!50.1960784313725!black, dotted]
table {%
1 1.39
2 2.78
4 5.56
8 11.12
16 22.24
};
\addlegendentry{Ideal quadratic complexity scaling}
\draw (axis cs:1,1.39) ++(0pt,2pt) node[
  scale=0.5,
  anchor=base west,
  text=black,
  rotate=0.0
]{1.4};
\draw (axis cs:2,1.68) ++(0pt,2pt) node[
  scale=0.5,
  anchor=base west,
  text=black,
  rotate=0.0
]{1.7};
\draw (axis cs:4,3.65) ++(0pt,2pt) node[
  scale=0.5,
  anchor=base west,
  text=black,
  rotate=0.0
]{3.6};
\draw (axis cs:8,4.71) ++(0pt,2pt) node[
  scale=0.5,
  anchor=base west,
  text=black,
  rotate=0.0
]{4.7};
\draw (axis cs:16,5.77) ++(0pt,2pt) node[
  scale=0.5,
  anchor=base west,
  text=black,
  rotate=0.0
]{5.8};
\end{axis}

\end{tikzpicture}}
\caption{$\mathsf{GOFMM}$ Compression }
    \label{fig:weaka}
	\end{subfigure}
		\begin{subfigure}{.5\textwidth}
\resizebox{0.98\textwidth}{!}{
\begin{tikzpicture}

\definecolor{color0}{rgb}{0.12156862745098,0.466666666666667,0.705882352941177}

\begin{axis}[tuftelike,
legend cell align={left},
legend style={
  fill opacity=0.8,
  draw opacity=1,
  text opacity=1,
  at={(0.03,0.97)},
  anchor=north west,
  draw=white!80!black
},
log basis x={10},
tick align=outside,
tick pos=both,
x grid style={white!69.0196078431373!black},
xlabel={Nodes},
xmin=0.870550563296124, xmax=18.3791736799526,
xmode=log,
xtick style={color=black},
xtick={1,2,4,8,16},
xticklabels={1,2,4,8,16}, 
y grid style={white!69.0196078431373!black},
ylabel={Runtime (s)},
ymin=0.025, ymax=1.675,
ytick style={color=black},
axis y line*=left,  
]
\addplot [semithick, color0, mark=x, mark size=3, mark options={solid}]
table {%
1 0.1
2 0.13
4 0.19
8 0.23
16 0.26
};
\addlegendentry{$\mathsf{GOFMM}$ evaluation}
\addplot [semithick, white!50.1960784313725!black, dashed]
table {%
1 0.1
2 0.1
4 0.1
8 0.1
16 0.1
};
\addlegendentry{Ideal linear complexity scaling}
\addplot [semithick, white!50.1960784313725!black, dotted]
table {%
1 0.1
2 0.2
4 0.4
8 0.8
16 1.6
};
\addlegendentry{Ideal quadratic complexity scaling}
\draw (axis cs:1,0.1) ++(0pt,2pt) node[
  scale=0.5,
  anchor=base west,
  text=black,
  rotate=0.0
]{0.10};
\draw (axis cs:2,0.13) ++(0pt,2pt) node[
  scale=0.5,
  anchor=base west,
  text=black,
  rotate=0.0
]{0.13};
\draw (axis cs:4,0.19) ++(0pt,2pt) node[
  scale=0.5,
  anchor=base west,
  text=black,
  rotate=0.0
]{0.19};
\draw (axis cs:8,0.23) ++(0pt,2pt) node[
  scale=0.5,
  anchor=base west,
  text=black,
  rotate=0.0
]{0.23};
\draw (axis cs:14,0.26) ++(0pt,2pt) node[
  scale=0.5,
  anchor=base west,
  text=black,
  rotate=0.0
]{0.26};
\end{axis}

\end{tikzpicture}}
\caption{$\mathsf{GOFMM}$ Evaluation }
\label{fig:weakb}
	\end{subfigure}
	\caption{Weak scaling measurements of a Gaussian kernel matrices generated
synthetically with 6-D point clouds with roughly 6.25k$\times$6.25k with 1 node, 12.5k$\times$12.5k with 2 nodes, up to {$\scriptstyle\sim$}100k$\times$100k with 16 nodes.  Memory and runtime for exact multiplication of a dense matrix scales quadratically, hence the 2-node problem would correspond to 4-times the memory/computational cost (dotted) in total. Note that the above figures have log scale on the x-axis and linear scale on the y-axis. 
}
\label{fig:weak}
	\end{figure}

In weak scaling, the computational effort per resource (core/node) stays constant. We scale the problem size with respect to nodes and hence, for algorithms with linear complexity, the problem size per node stays constant.
But since matrix size scales quadratically, doubling the problem size would require that we scale the number of nodes quadratically in order to maintain a constant computational load per node. This quickly becomes infeasible due to limited computational resources. Therefore we scale the nodes linearly instead and provide corresponding ideal runtimes through the dotted lines in \autoref{fig:weak}.
  .\footnote{In theory, for a matrix with an off-diagonal rank of \(r_O\) , $\mathsf{GOFMM}$ has a computational complexity of \(\mathcal{O}(N\cdot r_O)\). But with certain adaptive rank selection and a certain accuracy, it potentially increases with problem size and thus for simplicity, we refrain to  \(\mathcal{O}(N \log N)\). }

 Weak scaling for $\mathsf{GOFMM}$ compression in \autoref{fig:weaka} results in a runtime complexity between \(\mathcal{O}(N)\) (linear) and \(\mathcal{O}(N^2)\) (quadratic). Although with this inaccurate complexity estimate we cannot measure the parallel efficiency and communication overhead of $\mathsf{GOFMM}$, it still shows us that it scales really well with increasing problem size and thereby proving that \(\mathcal{H}\)-matrix approximation is very \textit{beneficial} for large matrices compared to an exact dense multiplication.  \autoref{fig:weakb} shows runtimes for matrix-multiplication (also referred to as evaluation) with $\mathsf{GOFMM}$ for increasing problem size and nodes. We observe that the runtime for a problem size of $6.25k$ with 1 node is $0.10s$ and for a problem size of\(6.25k*16\approx100k\) with 16 nodes is about $0.26s$. Assuming \(\mathcal{O}(N \log N)\) computational complexity, ideal scaling would result in a runtime of ${0.10s} *\log (16) \approx {0.12s}$. Instead, the runtime of $0.26s$ we obtained results in a parallel efficiency of $\frac{0.12}{0.26}\approx50\%$.


 To summarize, we see a difference between $\mathsf{GOFMM}$'s \(\mathcal{O}(N \log N)\) runtime complexity and a quadratic complexity for  large matrices with sizes above $25k\times25k$ (see behavior in \autoref{fig:weak}).

\subsubsection{Strong scaling}
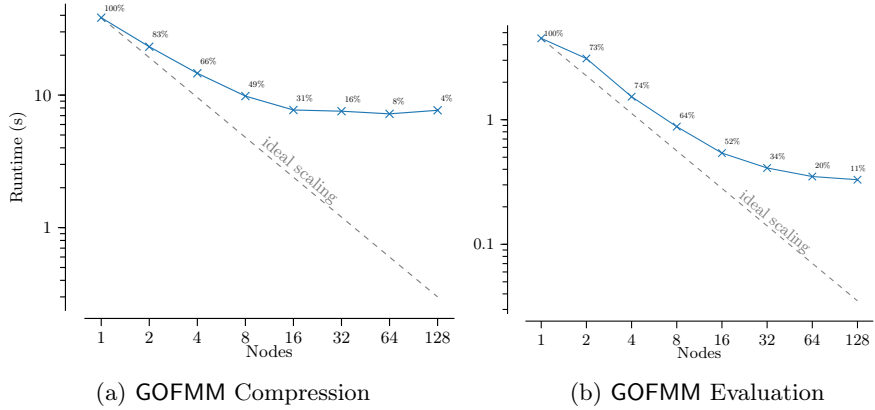
\begin{figure}
	\begin{subfigure}{.5\textwidth}
\resizebox{\textwidth}{!}{
\begin{tikzpicture}
\definecolor{color0}{rgb}{0.12156862745098,0.466666666666667,0.705882352941177}
\begin{axis}[tuftelike,
legend cell align={left},
legend style={fill opacity=0.8, draw opacity=1, text opacity=1, draw=white!80!black},
log basis x={10},
log basis y={10},
tick align=outside,
tick pos=both,
x grid style={white!69.0196078431373!black},
xlabel={Nodes},
xmin=0.784584097896751, xmax=163.143760296865,
xmode=log,
xtick style={color=black},
xtick={1,2,4,8,16,32,64,128},
xticklabels={1,2,4,8,16,32,64,128}, 
y grid style={white!69.0196078431373!black},
ylabel={Runtime (s)},
ymin=0.235375229369025, ymax=48.9431280890596,
ymode=log,
ytick style={color=black},
            yticklabels={0.1,1,10,20,30},
axis y line*=left,  
]
\addplot [semithick, color0, mark=x, mark size=3, mark options={solid}]
table {%
1 38.4
2 23.18
4 14.62
8 9.83
16 7.72
32 7.56
64 7.22
128 7.68
};
\addplot [semithick, white!50.1960784313725!black, dashed]
table {%
1 38.4
2 19.2
4 9.6
8 4.8
16 2.4
32 1.2
64 0.6
128 0.3
};
\node[rotate=-40, text=gray] at (axis cs:17.5,2.7) {ideal scaling};

\draw (axis cs:1,36.5) ++(0pt,5pt) node[
  scale=0.5,
  anchor=base west,
  text=black,
  rotate=0.0
]{100\%};
\draw (axis cs:2,23.18) ++(0pt,5pt) node[
  scale=0.5,
  anchor=base west,
  text=black,
  rotate=0.0
]{83\%};
\draw (axis cs:4,14.62) ++(0pt,5pt) node[
  scale=0.5,
  anchor=base west,
  text=black,
  rotate=0.0
]{66\%};
\draw (axis cs:8,9.83) ++(0pt,5pt) node[
  scale=0.5,
  anchor=base west,
  text=black,
  rotate=0.0
]{49\%};
\draw (axis cs:16,7.72) ++(0pt,5pt) node[
  scale=0.5,
  anchor=base west,
  text=black,
  rotate=0.0
]{31\%};
\draw (axis cs:32,7.56) ++(0pt,5pt) node[
  scale=0.5,
  anchor=base west,
  text=black,
  rotate=0.0
]{16\%};
\draw (axis cs:64,7.22) ++(0pt,5pt) node[
  scale=0.5,
  anchor=base west,
  text=black,
  rotate=0.0
]{8\%};
\draw (axis cs:128,7.68) ++(0pt,5pt) node[
  scale=0.5,
  anchor=base west,
  text=black,
  rotate=0.0
]{4\%};
\end{axis}

\end{tikzpicture}}
\caption{$\mathsf{GOFMM}$ Compression }
   \label{fig:stronga}
	\end{subfigure}
	\begin{subfigure}{.5\textwidth}
	\resizebox{0.9\textwidth}{!}{
\begin{tikzpicture}

\definecolor{color0}{rgb}{0.12156862745098,0.466666666666667,0.705882352941177}

\begin{axis}[tuftelike,
legend cell align={left},
legend style={fill opacity=0.8, draw opacity=1, text opacity=1, draw=white!80!black},
log basis x={10},
log basis y={10},
tick align=outside,
tick pos=both,
x grid style={white!69.0196078431373!black},
xlabel={Nodes},
xmin=0.784584097896751, xmax=163.143760296865,
xmode=log,
xtick style={color=black},
xtick={1,2,4,8,16,32,64,128},
xticklabels={1,2,4,8,16,32,64,128}, 
            yticklabels={0.01,0.1,1,2,3,4,5},
y grid style={white!69.0196078431373!black},
ymin=0.0275830346916826, ymax=5.73552282293668,
ymode=log,
ytick style={color=black},
axis y line*=left,  
]
\addplot [semithick, color0, mark=x, mark size=3, mark options={solid}]
table {%
1 4.5
2 3.1
4 1.53
8 0.88
16 0.54
32 0.41
64 0.35
128 0.33
};
\addplot [semithick, white!50.1960784313725!black, dashed]
table {%
1 4.5
2 2.25
4 1.125
8 0.5625
16 0.28125
32 0.140625
64 0.0703125
128 0.03515625
};
\node[rotate=-40, text=gray] at (axis cs:36,0.15) {ideal scaling};
\draw (axis cs:1,4) ++(0pt,5pt) node[
  scale=0.5,
  anchor=base west,
  text=black,
  rotate=0.0
]{100\%};
\draw (axis cs:2,3.1) ++(0pt,5pt) node[
  scale=0.5,
  anchor=base west,
  text=black,
  rotate=0.0
]{73\%};
\draw (axis cs:4,1.53) ++(0pt,5pt) node[
  scale=0.5,
  anchor=base west,
  text=black,
  rotate=0.0
]{74\%};
\draw (axis cs:8,0.88) ++(0pt,5pt) node[
  scale=0.5,
  anchor=base west,
  text=black,
  rotate=0.0
]{64\%};
\draw (axis cs:16,0.54) ++(0pt,5pt) node[
  scale=0.5,
  anchor=base west,
  text=black,
  rotate=0.0
]{52\%};
\draw (axis cs:32,0.41) ++(0pt,5pt) node[
  scale=0.5,
  anchor=base west,
  text=black,
  rotate=0.0
]{34\%};
\draw (axis cs:64,0.35) ++(0pt,5pt) node[
  scale=0.5,
  anchor=base west,
  text=black,
  rotate=0.0
]{20\%};
\draw (axis cs:110,0.33) ++(0pt,5pt) node[
  scale=0.5,
  anchor=base west,
  text=black,
  rotate=0.0
]{11\%};
\end{axis}

\end{tikzpicture}}
 \caption{$\mathsf{GOFMM}$ Evaluation }
 \label{fig:strongb}
	\end{subfigure}
\caption{Strong scaling measurements of a Gaussian kernel matrices generated
synthetically with 6-D point clouds, all roughly of size 100k-by-100k. Next to the data cross is the parallel efficiency in percent. Results run on Skylake partition of SuperMuc-NG. Each node has 48 cores, 128 Nodes hence corresponds to 6144 cores. }
\label{fig:gofmm_strong}
 \end{figure}

In strong scaling, the problem size stays constant while increasing the computational resources and this can be challenging due to diminishing computational work per node and increasing communication overhead.
\\
\\
 \autoref{fig:gofmm_strong} shows strong scaling measurements for $\mathsf{GOFMM}$ compression and evaluation for a $100k\times 100k$ synthetic kernel matrix. 
 \autoref{fig:stronga} on the \textbf{left} we see the one-time compression time (For parameter see\footnote{$\mathsf{GOFMM}$ parameters: $\texttt{max\_leaf\_node\_size}=768$, $\texttt{max\_off\_diagonal\_ranks}=768$, $\texttt{user\_tolerance}=1E-3$, $\texttt{num\_neighbors}=64$}). 
 Compression algorithm for a 6D random Gaussian kernel matrix of size $100k\times 100k$ takes \textit{13s} on one node while multiplication with a vector of size $100k\times 512$ has a runtime of \textit{1.35s}. We can also observe that the parallel efficiency for both algorithms ranges down to 4\% and 11\% with 128 nodes and that there is no performance gain when nodes higher than 16 are used. As mentioned previously, it is not unusual for efficiency to have tendencies of stagnation or deterioration with strong scaling due to problems such as increasing communication overhead and load imbalance.

Having a limit on maximum acceptable efficiency is not unusual for parallel code; also to reiterate, growth in runtimes are possible as communication times are increasing. For this reason we highlight similar runtime scaling for matrix evaluation and the one-time matrix compression cost also mentioned in \cite{mpigofmm}.

We see a similar tendency in \autoref{fig:strongb} starting with 52\% efficiency with 16 nodes, implying that the runtime is 8-times slower than on a single node. 
Note that we also run a problem size of $100k\times 100k$ with 16 nodes for weak scaling in \autoref{fig:weak} and get similar results as expected. 

\section{Conclusion}
With ever-growing applications with non-linear high-dimensional data in Machine learning and AI, it becomes more and more difficult to process this data efficiently. We utilize a manifold learning algorithm (of $\mathsf{datafold}$) to compute the underlying lower dimension of such data and propose an approach to reduce the computational complexity of certain operations contained in such algorithms. We present a proof-of-concept that hierarchical methods can be applied to large matrices in aforementioned algorithms. Since $\mathsf{datafold}$ is written in Python and $\mathsf{GOFMM}$ is written in C++, the overhead caused by the $\mathsf{SWIG}$ interface are unknown. This also causes more limitations on the ability to fully utilize $\mathsf{GOFMM}$'s MPI functionality. In ongoing work we integrate other kernels that require dense matrices, and thus are more suitable to the approach and make use of $\mathsf{GOFMM}$ to its full potential. 


\subsubsection{Acknowledgements} 

This  material  is  based  upon  work  supported  by the Technical University of Munich, by the International Graduate School of Science and Engineering (IGSSE) with funds from Deutsche Forschungsgemeinschaft (DFG) through SPPEXA~\cite{bungartz2020software} and by the Competence Network for Scientific High Performance Computing in Bavaria (KONWIHR) with funds from  Bayerisches Staatsministerium für Wissenschaft und Kunst (STMWK).

%
%
%
%


\begin{thebibliography}{10}
	\providecommand{\url}[1]{\texttt{#1}}
	\providecommand{\urlprefix}{URL }
	\providecommand{\doi}[1]{https://doi.org/#1}
	
	\bibitem{17}
	ARPACK Software. http://www.caam.rice.edu/software/ARPACK/
	
	\bibitem{10}
	Bebendorf, M.: Hierarchical matrices. Springer Publishing Company,
	Incorporated, 1st edition (2008)
	
	\bibitem{bolager2023sampling}
	Bolager, E.L., Burak, I., Datar, C., Sun, Q., Dietrich, F.: Sampling weights of
	deep neural networks. arXiv preprint arXiv:2306.16830  (2023)
	
	\bibitem{bungartz2020software}
	Bungartz, H.J., Nagel, W.E., Neumann, P., Reiz, S., Uekermann, B.: {Software
		for Exascale Computing: Some Remarks on the Priority Program SPPEXA}. In:
	Bungartz, H.J., Reiz, S., Uekermann, B., Neumann, P., Nagel, W.E. (eds.)
	Software for Exascale Computing - SPPEXA 2016-2019. pp. 3--18. Springer
	International Publishing, Cham (2020)
	
	\bibitem{COIFMAN20065}
	Coifman, R.R., Lafon, S.: Diffusion maps. Applied and Computational Harmonic
	Analysis  \textbf{21},  5--30 (2006).
	\doi{https://doi.org/10.1016/j.acha.2006.04.006},
	\url{https://www.sciencedirect.com/science/article/pii/S1063520306000546},
	special Issue: Diffusion Maps and Wavelets
	
	\bibitem{deng2012mnist}
	Deng, L.: The mnist database of handwritten digit images for machine learning
	research. IEEE Signal Processing Magazine  \textbf{29}(6),  141--142 (2012)
	
	\bibitem{hessianembed}
	Donoho, D.L., Grimes, C.: Hessian eigenmaps: Locally linear embedding
	techniques for high-dimensional data. Proceedings of the National Academy of
	Sciences  \textbf{100}(10),  5591--5596 (2003).
	\doi{10.1073/pnas.1031596100},
	\url{https://www.pnas.org/doi/abs/10.1073/pnas.1031596100}
	
	\bibitem{garriga2018deep}
	Garriga-Alonso, A., Rasmussen, C.E., Aitchison, L.: Deep convolutional networks
	as shallow gaussian processes. arXiv preprint arXiv:1808.05587  (2018)
	
	\bibitem{grasedyck2013literature}
	Grasedyck, L., Kressner, D., Tobler, C.: A literature survey of low-rank tensor
	approximation techniques. GAMM-Mitteilungen  \textbf{36},  53--78 (2013)
	
	\bibitem{11}
	Hackbusch, W.: Hierarchical Matrices: Algorithms and Analysis. Springer-Verlag
	Berlin Heidelberg (2015)
	
	\bibitem{doi:10.1137/100804139}
	Halko, N., Martinsson, P.G., Shkolnisky, Y., Tygert, M.: An algorithm for the
	principal component analysis of large data sets. SIAM Journal on Scientific
	Computing  \textbf{33}(5),  2580--2594 (2011). \doi{10.1137/100804139},
	\url{https://doi.org/10.1137/100804139}
	
	\bibitem{3}
	Hofmann, T., Schölkopf, B., Smola, A.J.: Kernel methods in machine learning.
	The annals of statistics (2008), 1171–1220 (2008)
	
	\bibitem{5}
	Lehmberg, D., Dietrich, F., Köster, G., Bungartz, H.J.: datafold: data-driven
	models for point clouds and time series on manifolds. Journal of Open Source
	Software, 5(51), 2283 (2020)
	
	\bibitem{19}
	Lehoucq, R.B., Sorensen, D.C., Yang, C.: ARPACK USERS GUIDE: Solution of Large
	Scale Eigenvalue Problems by Implicitly Restarted Arnoldi Methods. SIAM,
	Philadelphia, PA (1998)
	
	\bibitem{doi:10.1137/20M1349667}
	Liu, Y., Ghysels, P., Claus, L., Li, X.S.: Sparse approximate multifrontal
	factorization with butterfly compression for high-frequency wave equations.
	SIAM Journal on Scientific Computing  \textbf{43}(5),  S367--S391 (2021).
	\doi{10.1137/20M1349667}, \url{https://doi.org/10.1137/20M1349667}
	
	\bibitem{MARTINSSON201147}
	Martinsson, P.G., Rokhlin, V., Tygert, M.: A randomized algorithm for the
	decomposition of matrices. Applied and Computational Harmonic Analysis
	\textbf{30}(1),  47--68 (2011).
	\doi{https://doi.org/10.1016/j.acha.2010.02.003},
	\url{https://www.sciencedirect.com/science/article/pii/S1063520310000242}
	
	\bibitem{scikit-learn}
	Pedregosa, F., Varoquaux, G., Gramfort, A., Michel, V., Thirion, B., Grisel,
	O., Blondel, M., Prettenhofer, P., Weiss, R., Dubourg, V., Vanderplas, J.,
	Passos, A., Cournapeau, D., Brucher, M., Perrot, M., Duchesnay, E.:
	Scikit-learn: Machine learning in {P}ython. Journal of Machine Learning
	Research  \textbf{12},  2825--2830 (2011)
	
	\bibitem{reiz2022neural}
	Reiz, S., Neckel, T., Bungartz, H.J.: Neural nets with a newton conjugate
	gradient method on multiple gpus. In: International Conference on Parallel
	Processing and Applied Mathematics. pp. 139--152. Springer (2022)
	
	\bibitem{locallylinearembed}
	Roweis, S.T., Saul, L.K.: Nonlinear dimensionality reduction by locally linear
	embedding. Science  \textbf{290}(5500),  2323--2326 (2000).
	\doi{10.1126/science.290.5500.2323},
	\url{https://www.science.org/doi/abs/10.1126/science.290.5500.2323}
	
	\bibitem{18}
	Sorensen, D.C.: Implicitly restarted Arnoldi/Lanczos methods for large scale
	eigenvalue calculations. SIAM J. Matrix Anal. Appl., 13, pp. 357–385 (1992)
	
	\bibitem{isomap}
	Tenenbaum, J., Silva, V., Langford, J.: A global geometric framework for
	nonlinear dimensionality reduction. Science  \textbf{290},  2319--2323 (01
	2000)
	
	\bibitem{16}
	Virtanen, P., Gommers, R., Oliphant, T.E., Haberland, M., Reddy, T., et~al.:
	{{SciPy} 1.0: Fundamental Algorithms for Scientific Computing in Python}.
	Nature Methods (2020)
	
	\bibitem{1}
	Yu, C.D., Levitt, J., Reiz, S., Biros, G.: Geometry- Oblivious FMM for
	Compressing Dense SPD Matrices. In Proceedings of SC17, Denver, CO, USA
	(2017)
	
	\bibitem{mpigofmm}
	Yu, C.D., Reiz, S., Biros, G.: Distributed-memory hierarchical compression of
	dense spd matrices. In: SC18: International Conference for High Performance
	Computing, Networking, Storage and Analysis. pp. 183--197 (2018).
	\doi{10.1109/SC.2018.00018}
	
	\bibitem{23}
	Yu, C.D., Reiz, S., Biros, G.: Distributed O(N) Linear Solver for Dense
	Symmetric Hierarchical Semi-Separable Matrices. IEEE 13th International
	Symposium on Embedded Multicore/Many-core Systems-on-Chip (MCSoC) (2019)
	
\end{thebibliography}

\end{document}